\documentclass[12pt]{amsart}
\usepackage{enumerate}

\newcommand{\E}{\mathcal E}
\newcommand{\R}{\mathbb R}
\newcommand{\C}{\mathbb C}
\renewcommand{\th}{\theta}
\newcommand{\Th}{\Theta}
\newcommand{\g}{\mathfrak g}

\newcommand{\sll}{\mathfrak{sl}}
\newcommand{\gl}{\mathfrak{gl}}
\newcommand{\RP}{\R\mathrm{P}}
\newcommand{\CP}{\C\mathrm{P}}
\newcommand{\W}{\mathcal W}
\newcommand{\p}{\partial}
\newcommand{\cF}{\mathcal F}
\newcommand{\eps}{\epsilon}
\renewcommand{\O}{\mathcal O}
\renewcommand{\L}{\mathcal L}
\renewcommand{\P}{\mathrm{P}}
\renewcommand{\S}{\mathcal S}
\newcommand{\tr}{\operatorname{tr}}
\newcommand{\End}{\operatorname{End}}
\newcommand{\Gr}{\operatorname{Gr}}
\newcommand{\pp}[1]{\frac{\partial}{\partial #1}}
\newcommand{\om}{\omega}
\newcommand{\ad}{\operatorname{ad}}
\renewcommand{\Im}{\operatorname{Im}}

\newtheorem{thm}{Theorem}

\theoremstyle{definition}
\newtheorem{ex}{Example}

\theoremstyle{remark}
\newtheorem{rem}{Remark}

\title[Generalized Wilczynski invariants for non-linear ODEs]{Generalized Wilczynski
invariants for non-linear ordinary differential equations}
\author{Boris Doubrov}
\address{Belarussian State University, Skoriny 4, 220050, Minsk, Belarus}
\email{doubrov@islc.org} 

\begin{document}

\begin{abstract}
We show that classical Wilczynski--Se-ashi invariants of linear
systems of ordinary differential equations are generalized in a
natural way to contact invariants of non-linear ODEs. We explore
geometric structures associated with equations that have vanishing
generalized Wilczynski invariants and establish relationship of such
equations with deformation theory of rational curves on complex
algebraic surfaces.
\end{abstract}

\maketitle

\section{Introduction}
This paper is devoted to a very important class of contact invariants
of (non-linear) systems of ordinary differential equations. They can
be considered as a direct generalization of classical Wilczynski
invariants of linear differential equations, which are, in its turn,
closely related to projective invariants of non-parametrized curves.

The construction of these invariants is based on the fundamental idea
of approximating the non-linear objects by linear. In case of
differential equations we can consider linearization of a given
non-linear equation along each solution. Roughly speaking, this linearization
describes the tangent space to the solution set of the given
equation. Note that in general this set is not Hausdorff, but we can
still speak about linearization of the equation without any loss of
generality.

Unlike the class of all non-linear equations of a fixed order, which
is stable with respect contact transformations, the class of linear
equations forms a category with a much smaller set of morphisms. This
makes it possible to describe the set of all invariants of the linear
equations explicitly. The main goal of this paper is to show the
general geometric procedure that extends these invariants to the class
of non-linear equations via the notion of the linearization.

In short, the major result of this paper is that \emph{invariants of
the linearization of a given non-linear equation are contact invariants
of this equation}.

In fact, these invariants have been known for ordinary differential equations
of low order: third order equations~\cite{chern,sat:yos}, fourth
order ODEs~\cite{bry4ord,fels2,dun}, systems of second order~\cite{fels, gros}, but it
was not understood up to now that these invariants come from the linearization of
non-linear equations.

In the theory of linear differential equations Wilczynski invariants generate
all invariants and, in particular, are responsible for the trivialization
of the equation. In non-linear case the generalized Wilczynski invariants form
only part of the generators in the algebra of all contact invariants. However,
the equations with vanishing generalized Wilczynski invariants have a remarkable
property, that its solution space carries a natural geometric structure
(see~\cite{bry4ord} and \cite{dun}). We discuss this structure in more detail
in Section~\ref{sec:van} for a single ODE. Finally, in Section~\ref{sec:sys}
we show that all results of this paper can be extended to the
case of systems of ordinary differential equations.

\textbf{Acknowledgments.} This paper is based in large on the lectures given at
the IMA workshop ``Symmetries and overdetermined systems of partial
differential equations''. I would like to thank organizers of this workshop for
inviting me to give these lectures. I also would like to thank Eugene
Ferapontov, Rod Gover and Igor Zelenko for corrections and valuable discussions
on the topic of this paper.

\section{Naive approach}
As an example, let us outline this idea for the case of a
single non-linear ODE of order~$\ge3$. Note that all equations
of smaller order are contact equivalent to each other and, thus, do
not have any non-trivial invariants.

The complete set of invariants for a single linear ODE was described in
the classical work of Wilczynski~\cite{wil}. Namely, consider the class of linear
homogeneous differential equations
\begin{equation}\label{lin:eq}
y^{(n+1)}+p_n(x)y^{(n)}+\dots+p_0(x)y=0,
\end{equation}
viewed up to all invertible transformations of the form:
\begin{equation}\label{lin:transf}
(x,y)\mapsto(\lambda(x),\mu(x)y).
\end{equation}
In fact, these are the most general transformations preserving the
class of linear equations.

A function $I$ of the coefficients $p_i(x)$ and their
derivatives is called \emph{a relative invariant of weight $l$} if it is
transformed by the rule $I\mapsto(\lambda')^lI$ under the change
of variables~\eqref{lin:transf}. In particular, if such relative
invariant vanishes identically for the initial equation, it will also
vanish for the transformed one. Relative invariants of weight $0$ are
called \emph{(absolute) invariants of the linear equation}.

E.~Wilczynski~\cite{wil} gave the complete description of all relative
(and, thus, absolute) invariants of linear ODEs of any order. It is
well-known that each equation~\eqref{lin:eq} can be brought by the
transformation~\eqref{lin:transf} to the so-called Laguerre--Forsyth
canonical form:
\[
y^{(n+1)} + q_{n-2}(x)y^{(n-2)} + \dots + q_0(x)y=0.
\]
The set of transformations preserving this canonical form is already a
finite-dimensional Lie group:
\[
(x,y)\mapsto \left(\frac{ax+b}{cx+d}, \frac{ey}{(cx+d)^{n+1}}\right).
\]
This group acts on coefficients $q_0,\dots,q_{n-2}$ of the canonical
form, and the relative invariants of this action are identified with the
relative invariants of the general linear equation. The simplest $n-1$
relative invariants $\th_3,\dots,\th_{n+1}$ linear in
$q_i^{(j)}$ have the form:
\begin{equation}\label{eq:linv}
\th_k = \sum_{j=1}^{k-2} (-1)^{j+1}
\frac{(2k-j-1)!(n-k+j)!}{(k-j)!(j-1)!} q_{n-k+j}^{(j-1)} \qquad
k=3,\dots, n+1.
\end{equation}
Each invariant $\th_k$ has weight $k$. Wilczynski proved that all
other invariants can be expressed in terms of the invariants $\th_k$
and their derivatives.

Although these formulas express invariants in terms of coefficients of
the canonical form, they can be written explicitly in terms of the
initial coefficients of the general equation. Moreover, it can be
shown that each of these invariants $\th_i$ is polynomial in terms of
functions $p_i(x)$ and their derivatives.

Equivalence theory of linear ODEs is intimately related with the
projective theory of non-parametrized curves. Namely,
let $\{y_0(x),\dots,y_n(x)\}$ be a fundamental set of solutions of
a linear equation~$\E$ given by~\eqref{lin:eq}. Consider the curve
\[
L_{\E}=\{[y_0(x):y_1(x):\dots:y_n(x)]\mid x\in \R\}
\]
in the $n$-dimensional projective space
$\RP^n$. Since the solutions do not vanish simultaneously for any
$x\in\R$, this curve is well-defined. Moreover, since the fundamental
set of solutions is linearly independent, this curve is
\emph{non-degenerate}, i.e., it is not contained in any hyperplane.
Finally, since the set of fundamental solutions is defined up to any
non-degenerate linear transformation, we see that the curve $L_{\E}$
is defined up to projective transformations.

Consider now what happens with this curve, if we apply the
transformations~\eqref{lin:transf} to the equation~$\E$. The
transformations $(x,y)\mapsto(x,\mu(x)y)$ do not change the curve,
since they just multiply each solution by $\mu(x)$. The
transformations $(x,y)\mapsto(\lambda(x),y)$ are equivalent to
reparametrizations of $L_{\E}$.

Thus, we see that to each linear equation $\E$ we can assign the set
of projectively-equivalent non-degenerate curves in $\RP^n$. It is
easy to see that this correspondence is one-to-one. Indeed, having a
non-degenerate curve in $\RP^n$ we can fix a parameter $x$ on it and
write it explicitly as $[y_0(x):y_1(x):\dots:y_n(x)]$, where the
coordinates $y_i(x)$ are defined modulo a non-zero multiplier
$\mu(x)$. Since the curve is non-degenerate, the set of functions
$\{y_i(x)\}_{i=0,\dots,n}$ is linearly independent and defines a
unique linear equation $\E$ having these functions as the set of
fundamental solutions.

Each relative invariant $I$ of weight $k$ can be naturally interpreted
as a section of the line bundle $S^k(TL)$ invariant with respect to
projective transformations. In particular, relative invariants of
weight $0$ are just projective differential invariant of
non-parametrized curves in the projective space. First examples of
such invariants were constructed by Sophus Lie~\cite{lie} for curves
on the projective plane and then generalized by Halphen~\cite{halph}
to the case of projective spaces of higher dimensions. Note that they
can be constructed via the standard Cartan moving frame method.

E.~Wilczynski also proved the following result characterizing the
equations with vanishing invariants:
\begin{thm}\label{wil:thm}
  Let $\E$ be a linear ODE given by~\eqref{lin:eq}. The following
  conditions are equivalent:
  \begin{enumerate}
  \item invariants $\th_3,\dots,\th_{n+1}$ vanish identically.
  \item the equation $\E$ is equivalent to the trivial equation
    $y^{(n+1)}=0$.
  \item the curve $L_{\E}$ is an open part of the normal rational
    curve in $\RP^n$.
  \item the symmetry algebra of $L_{\E}$ is isomorphic to the
    subalgebra $\sll(2,\R)\subset\sll(n+1,\R)$ acting irreducibly on
    $\R^{n+1}$.
  \end{enumerate}
\end{thm}

Let us show how to extend Wilczynski invariants to arbitrary
non-linear ordinary differential equations via the notion of
linearization. Indeed, consider now an arbitrary non-linear ODE solved
with respect to the highest derivative:
\begin{equation}\label{eq:nonlin}
y^{(n+1)} = f(x,y,y',\dots,y^{(n)}).
\end{equation}

Let $y_0(x)$ be any solution of this equation. Then we can consider
the linearization of~\eqref{eq:nonlin} along this solution. It is a
linear equation
\begin{equation}\label{eq:liner}
  h^{(n+1)}=\frac{\p f}{\p y^{(n)}} h^{(n)}+\dots+\frac{\p f}{\p y}h,
\end{equation}
where all coefficients are evaluated at the solution $y_0(x)$. It
describes all deformations $y_\eps(x)=y_0(x)+\eps h(x)$ of the solution
$y_0(x)$, which satisfy the equation~\eqref{eq:nonlin} modulo
$o(\eps)$.

Consider now Wilczynski invariants $\th_3,\dots,\th_{n+1}$ of the
linearization~\eqref{eq:liner}. They are polynomial in terms of the
coefficients $\frac{\p f}{\p y^{(i)}}$ and their derivatives (evaluated
at the solution $y_0(x)$). In general, let $F(x,y,y',\dots,y^{(n)})$ be
any function of $x$, $y(x)$ and its derivatives evaluated at any
solution of the equation~\eqref{eq:nonlin}. Then, differentiating it
by $x$ is equivalent to applying the operator of total derivative:
\[
D=\pp{x} + y'{\pp y} + \dots + y^{(n)}\pp{y^{(n-1)}}+f\pp{y^{(n)}}.
\]
Thus, analytically, Wilczynski invariants of~\eqref{eq:liner} can be
expressed as polynomials in terms of $D^j\left(\frac{\p f}{\p
    y^{(i)}}\right)$ ($0\le i\le n$, $j\ge0$) and are independent of the
solution $y_0(x)$ we started with. More precisely, there are $(n-1)$
well-defined expressions $W_3$, \dots, $W_{n+1}$ polynomial in terms
of $D^j\left(\frac{\p f}{\p y^{(i)}}\right)$, which, being evaluated
at any solution $y_0(x)$ of the equation~\eqref{eq:nonlin}, give
Wilczynski invariants of its linearization along $y_0(x)$. We call
these expressions $W_3,\dots,W_{n+1}$ the \emph{generalized Wilczynski
  invariants of non-linear ordinary differential equations}.

To understand, what geometric objects correspond to generalized
Wilczynski invariants we need to consider the jet interpretation of
differential equations. Let $J^i=J^i(\R^2)$ be the space of all jets of
order $i$ of (non-parametrized) curves on the plane. Then for
each curve $L$ on the plane we can define its lift $L^{(i)}$ to the
jet space $J^i$ of order $i$. Then the equation~\eqref{eq:nonlin} can
be considered as a submanifold $\E\subset J^{n+1}$ of codimension~1,
and, by the existence and uniqueness theorem for ODEs, all lifts of
its solutions form a one-dimensional foliation on $\E$. Let us denote
its tangent one-dimensional distribution by $E$. Let $(x,y)$ be any
local coordinate system on the plane. Then it naturally defines a
local coordinate system $(x,y_0,y_1,\dots,y_i)$ on $J^i$ such that the
lift of the graph $(x,y(x))$ is equal to
$(x,y(x),y'(x),\dots,y^{(i)})$. In these coordinates the submanifold
$\E$ is given by equation $y_{n+1}=f(x,y_0,y_1,\dots,y_n)$, functions
$(x,y_0,y_1,\dots,y_n)$ form a local coordinate system on $\E$ and the
distribution $E$ is generated by the vector field
\[
D=\pp{x} + y_1\pp{y_0} + \dots + y_n\pp{y_{n-1}}+f{\pp y_n}.
\]
The expressions $W_3,\dots,W_{n+1}$ can be naturally interpreted as
functions on the equation manifold $\E$. In fact, each invariant
$W_i$, being a relative invariant of weight $i$ of all linearizations,
defines a section $\W_i$ of the line bundle $S^iE^*$ by the formula:
\[
\W_i(D,D,\dots,D)=W_i.
\]

The largest set of invertible transformations preserving the class of
ODEs of fixed order $i$ consists of so-called contact transformations,
which are the most general transformations of $J^i$ preserving the
class of lifts of plane curves. The main result of this paper can be
formulated as follows:
\begin{thm}\label{thm:winv}
  The sections $\W_i$, $i=3,\dots,n+1$ of the line bundles $S^iE^*$ are
  invariant with respect to contact transformations of $J^{n+1}$.
\end{thm}

\section{Algebraic model of Wilczynski invariants}

In this section we give an alternative algebraic description of
Wilczynski invariants, which is due to Se-ashi~\cite{se1,se2}. There
are two main reasons for providing an algebraic picture behind
Wilczynski invariants. First, it can be easily generalized to
invariants of systems of ODEs or even to more general classes of
linear finite type equations (see~\cite{se1}). Wilczynski
himself described only invariants for a single ODE of arbitrary order
and for linear systems of second order ODEs. It seems that analytic
methods become too elaborate to proceed with systems of higher order.

Second, Se-ashi's construction gives alternative analytic formulas for computing
Wilczynski invariants, which are independent of Laguerre-Forsyth
canonical form. In particular, this explains why Wilczynski invariants
are polynomial in the initial coefficients and their derivatives,
while the coefficients of the canonical form are not.

Denote by $V_n$ the set $S^n(\R^2)$ of all homogeneous polynomials of
degree $n$ in two variables $v_1,v_2$. The standard $GL(2,\R)$-action on
$\R^2$ is naturally extended to $V_n$ and turns it into an irreducible
$GL(2,\R)$-module. Denote by $\rho_n\colon GL(2,\R)\to GL(V_n)$ the
corresponding representation mapping.

Consider the corresponding action of $\gl(2,\R)$ on $V_n=S^n(\R^2)$. Denote
by $X,Y,H,Z$ the following basis in $\gl(2,\R)$:
\[
X=\left(\begin{smallmatrix}0 & 1 \\ 0 & 0 \end{smallmatrix}\right),
H=\left(\begin{smallmatrix}1 & 0 \\ 0 & -1 \end{smallmatrix}\right),
Y=\left(\begin{smallmatrix}0 & 0 \\ 1 & 0 \end{smallmatrix}\right),
Z=\left(\begin{smallmatrix}1 & 0 \\ 0 & 1 \end{smallmatrix}\right).
\]
Then the action of these basis elements on $V_n$ is equivalent to the
action of the following vector fields on $\R_n[v_1,v_2]$:
\[
X=v_1\pp{v_2},\ H=v_1\pp{v_1}-v_2\pp{v_2},\ Y=v_2\pp{v_1},\
Z=v_1\pp{v_1}+v_2\pp{v_2}.
\]
In the sequel we shall identify $\gl(2,\R)$ with its image in
$\gl(V_n)$ defined by this action.

We define a gradation on $V_n$ such that the polynomials $E_0=v_1^n$,
$E_1=v_1^{n-1}v_2$, \dots, $E_n=v_2^n$ have degrees $-n-1$, $-n$, \dots, $-1$
respectively. Then the elements $X$, $H$, $Y$, $Z$ define operators of
degrees $-1$, $0$, $1$ and $0$ respectively. Denote also by
$V_n^{(i)}$ the set of all elements in $V_n$ of degree $\le i$. These
subspaces define a filtration on $V_n$:
\[
\{0\}=V_n^{(-n-2)}\subset V_n^{(-n-1)}\subset\dots\subset V_n^{(-1)}=V_n.
\]

Let $E$ be a one-dimensional vector bundle over a one-dimensional
manifold~$M$ with a local coordinate~$x$. Denote by $J^n(E)$ the
$n$-th order jet bundle of $E$, which is a $(n+1)$-dimensional vector
bundle over~$M$. Then any $(n+1)$-th order linear homogeneous ODE can
be considered as a connection on $J^n(E)$ such that all its solutions,
being lifted to $J^n(E)$, are horizontal. Let $\cF(J^n(E))$ be the
frame bundle of $J^n(E)$. Since $\dim V_n = n+1$, we can identify
$\cF(J^n(E))$ as a set of all isomorphisms $\phi_x\colon V_n\to J_x^n(E)$.
This turns $\cF(J^n(E))$ into a principle $GL(V_n)$-bundle over~$M$.

Denote by $\om$ the corresponding connection form on $\cF(J^n(E))$.
In brief, the main idea of the Se-ashi works is that Wilczynski invariants can be
interpreted in terms of a natural reduction of the $GL(V_n)$-bundle
$\cF(J^n(E))$ to some $G$-subbundle~$P$ characterized by the following
conditions:
\begin{enumerate}
\item $G$ is the image of the lower-triangular matrices in $GL(2,\R)$
  under the representation $\rho_n\colon GL(2,\R)\to GL(V_n)$;
\item $\om|_P$ takes values in the subspace $\langle
  X,H,Z,Y,Y^2,\dots,Y^n\rangle$.
\end{enumerate}

The form $\om|_P$ can be decomposed into the sum of $\om_{\gl}$ with values in
$\gl(2,\R)\subset \gl(V_n)$ and $\sum_{i=2}^n \om_iY^i$. Then $\om_{\gl}$ defines a flat
projective structure on the manifold $M$, while the forms $\om_i$ (or, more
precisely, their values on the vector field $\pp{x}$) coincide up to
the constant with Wilczynski invariants $\th_{i+1}$.

Let us describe this reduction in more detail. Denote by $G^{(0)}$ the
subgroup of $GL(V_n)$ consisting of all elements preserving the
filtration $V_n^{(i)}$ on $V_n$ introduced above.  These are exactly
all elements of $GL(V_n)$ represented by lower-triangular matrices in the
basis $\{E_0,\dots,E_n\}$. For $k\ge1$ denote by $\gl^{(k)}(V_n)$ the
following subalgebra in $\gl(V_n)$:
\[
\gl^{(k)}(V_n)=\{\phi\in \gl(V_n)\mid \phi(V_n^{(i)})\subset V_n^{(i+k)}\}.
\]
Let $GL^{(k)}(V_n)$ be the corresponding unipotent subgroup in
$GL(V_n)$.  Define the subgroups $G^{(k)}\subset G^{(0)}$ as the
products $GL^{(0)}(2,\R)GL^{(k)}(V_n)$ for each $k\ge0$, where
$GL^{(0)}(2,\R)$ is the intersection of $G^{(0)}$ with
$\rho_n(GL(2,\R))$.  Denote also by $W$ the subspace in $\gl(V_n)$
spanned by the endomorphisms $Y^2,\dots,Y^n$ (here $Y\subset\gl(2,\R)$
is identified with the corresponding element in $\gl(V_n)$).

The reduction~$P\subset\cF(J^n(E))$ is constructed via series of
reductions $P_{k+1}\subset P_{k}$, where $P_k$ is a principal
$G^{(k)}$-bundle characterized by the following conditions:
\begin{enumerate}[a)]
\item $P_0$ consists of all frames $\phi_x\colon V_n\to J_x^n(E)$,
  which map the filtration of $V_n$ into the filtration on each fiber
  $J_x^n(E)$;
\item for $k\ge1$ the form $\om_k=\om|_{P_k}$ takes values in the subspace
  $W_k=W+\gl^{(k-1)}(V_n)+\gl(2,\R)$.
\end{enumerate}
At the end of this procedure we arrive at the principal bundle $P=P_{n+1}$ with
the structure group $G=G_{n+1}$ and the 1-from $\om=\om_{n+1}$ with
values in $W+\gl(2,\R)$.

In fact, the second condition can be considered as a definition of
$P_k$ for $k\ge1$. Indeed, let $(x,y_0,y_1,\dots,y_n)$ be a local
coordinate system on $J^n(E)$, such that the $n$-jet of the section
$y(x)$ of $J(E)$ is given by $y_i(x)=y^{(i)}(x)$. Then the connection
on $J^n(E)$ corresponding to the linear homogeneous
equation~\eqref{lin:eq} is defined as an annihilator of the forms:
\begin{align*}
\th_i&=dy_i-y_{i+1}dx,\\
\th_n&=dy_n+\left(\sum_{i=0}^np_i(x)y_i\right)dx.
\end{align*}
If $s\colon x\mapsto (y_0(x),\dots,y_n(x))$ is any section of
$J^n(E)$, then the covariant derivative of $s$ along the vector field
$\partial_x$ has the form:
\[
\nabla_{\partial_x} s =
(y_0'(x)-y_1(x),y_1'(x)-y_2(x),\dots,y_n'(x)+\sum_{i=0}^np_i(x)y_i(x)).
\]

Let $s_i$, $i=0,\dots,n$, be the standard sections defined by
$y_j(x)=\delta_{ij}$. Then $s=\{s_0,\dots,s_n\}$ is a local
section of the frame bundle $P_0$. Let $\om$ be the
connection form. The pull-back $s^*\om$ can be written in this
coordinate system as:
\[
s^*\om = \begin{pmatrix}
0 & -1 & 0  & \dots & 0\\
0 & 0  & -1 & \dots & 0\\
\vdots & \vdots & \vdots & \ddots & \vdots \\
0 & 0 & 0 & \dots & -1 \\
p_0(x) & p_1(x) & p_2(x) & \dots & p_n(x)
\end{pmatrix}dx
\]
Then Se-ashi reduction theorem says that there exists such
gauge transformation $C\colon M\to G^{(0)}$ that
\[
C^{-1}\big(s^*\om\big)C + C^{-1}dC = \big(-X+\alpha H +\beta Z +
\gamma Y + \sum_{i=2}^n\bar\th_{i+1}Y^i\big)\,dx,
\]
and the set of such transformations forms a principal $G$-bundle. As shown in~\cite{se1},
the coefficients $\bar\th_{i+1}$ coincide with classical Wilczynski
invariants $\th_{i+1}$ defined by~\eqref{eq:linv} up to the
constant and some polynomial expression of invariants of lower weight:
\[
\bar\th_{i+1} = c\th_{i+1} + P_{i+1}(\th_3,\dots,\th_i),
\]
where $P_{i+1}$ is some fixed polynomial without free term. In particular,
Theorem~\ref{wil:thm} remains true if we substitute invariants $\th_i$ with their
modified versions $\bar\th_i$.

We shall not repeat the computations from Se-ashi work~\cite{se1},
just mentioning that it is based on the following simple technical
fact. Namely, for any $k\ge1$ consider the subspace $\gl_k(V_n)$ of
all operators of degree $k$ and the mapping $\ad_k(X)\colon
\gl_k(V_n)\to\gl_{k-1}(V_n)$, $A\mapsto [X,A]$. Then we have the
decomposition $\gl_k(V_n)=\langle Y^k \rangle \oplus \Im \ad_{k+1}X$
for all $k\ge1$, which allows to carry effectively the reduction from
$P_k$ to $P_{k+1}$. On each step such reduction involves only the
operation of solving linear equations with constant coefficients and
differentiation. In particular, this proves, that Wilczynski
invariants are polynomial in terms of the coefficients
$p_0(x),\dots,p_n(x)$ of the initial equation~\eqref{lin:eq}.

In general, we can not make further reductions without assumption that
some of these invariants do not vanish. And if all these invariants
vanish, then our equation is equivalent to the trivial equation
$y^{(n+1)}=0$.

\section{Generalization of Wilczynski invariants to non-linear ODEs}

Let $\E$ be now an arbitrary non-linear ODE of order $(n+1)$. As above,
we identify it with a submanifold of codimension~$1$ in the jet space
$J^{n+1}$. Then in local coordinates $(x,y_0,\dots,y_{n+1})$ the
equation $\E$ is given by $y_{n+1}=f(x,y_0,\dots,y_n)$, and
the functions $(x,y_0,\dots,y_n)$ form a local coordinate system
on~$\E$.

Let us recall that there is a canonical contact $C^k$ distribution defined
on each jet space $J^k(\R^2)$, which is generated by tangent lines to
all lifts of curves from the plane. In local coordinates it is
generated by two vector fields:
\[
C^k=\left\langle \pp{y_k},\ \pp x+y_1\pp{y_0}+\dots y_k\pp{y_{k-1}}\right\rangle.
\]
We also have natural projections $\pi_{k,l}\colon
J^k(\R^2)\to J^l(\R^2)$ for all $l<k$.

The contact distribution $C^{n+1}$ defines a line bundle on the
equation~$\E$ as follows $E_p=T_p\E\cap C^{n+1}_p$ for all
$p\in\E$. Integral curves of this bundle are precisely lifts of the
solutions of the given ODE. In local coordinates this line bundle is
generated by the vector field:
\[
D = \pp x + y_1\pp{y_0}+\dots+y_n\pp{y_{n-1}}+f\pp{y_n},
\]
which defines also the total derivative operator.

We are interested in the frames on the normal bundle to $E$, that is
the bundle $N(\E)=T\E/E$. We have a natural filtration of this bundle
defined by means of the projections $\pi_{n+1,i}\colon\E\to J^i(\R^2)$
for $i<n$. Namely, we define $N_i$ as the intersection of $\pi_{n+1,i}^*C^i$
with $T\E$ modulo the line bundle $E$ for all $i=0,\dots,n$. Then it
is easy to see that the sequence
\[
N=N_0\supset N_1 \supset \dots N_n \supset 0,
\]
is strictly decreasing and $\dim N_i=n+1-i$. In local coordinates we
have
\[
N_i=\left\langle \pp{y_i},\dots,\pp{y_n}\right\rangle + E.
\]
The vector field $D$ defines a first order operator on $N$, which is
compatible with this filtration, i.e. $D(N_i)\subset N_{i-1}$ for all
$i=1,\dots,n$.

As in the case of linear equations, we can define the frame bundle,
consisting of all maps $\phi_p\colon V_n\to N_p$ from the
$GL(2,R)$-module $V_n$ into $N_p$, which preserve the filtrations.
This defines a $G^{(0)}$-bundle $P_0$ over $\E$.

To proceed further with similar reductions, we need to have an
analog of the connection form. In general, we don't have any natural
connection on the normal bundle $N(\E)$. However, we can define the
covariant derivative along all vectors lying in the line bundle $E$
generated by $D$:
\[
  \nabla_{D}(X + E) = [D,X] + E,\quad\text {for any } X + E \in N(\E).
\]
This gives us a so-called partial connection with the
connection form $\om\colon E \to \gl(V_n)$. This connection form is
sufficient for our purposes. As in the case of linear equations,
we can form the similar set of reductions $P_i$ of the principle bundle $P_0$.
Each such bundle $P_i$, being restricted to any solution $L$, coincides with the
principle bundle constructed from the linearization of $\E$ along this
solution. As a result, we get the set of well-defined generalized
Wilczynski invariants $\W_i$, where each $\W_i$ is a global section of
the bundle $S^iE^*$.

Since the bundle $N(\E)$ and the connection form $\om$ are defined
totally in terms of the contact geometry of the jet space $J^{n+1}$
and the reduction procedure does not depend on any external data, we
see that generalized Wilczynski invariants are contact invariants of
the original equation. It can be formulated rigorously in the
following way:
\begin{thm} Assume that $n\ge2$.  Let $\phi\colon J^{n+1}\to J^{n+1}$ be a local contact
  transformation establishing the local equivalence of two equations
  $\E$ and $\overline \E$. Let $E$ ($\overline E$) be the line bundle
  on $\E$ (resp.~$\overline\E$) defining the solution foliation. Let
  $\W_i\in \Gamma(S^iE^*)$ (resp.~$\overline \W_i\in
  \Gamma(S^i\overline E^*)$ be generalized Wilczynski invariants of
  $\E$ (resp.~$\overline\E$). Then we have $\phi_*(E)=\overline E$ and
  $\phi^*\overline \W_i=\W_i$ for all $i=3,\dots,n+1$.
\end{thm}

This theorem was first proven in the work~\cite{d2001} via direct
computation of the canonical Cartan connection associated with any
ordinary differential equation of order $\ge4$ (see~\cite{tan89,dkm} for more
details).

\section{Equations with vanishing Wilczynski invariants}
\label{sec:van}
\subsection{Structure of the solution space}
Suppose now that all generalized Wilczynski invariants vanish identically for a
given equation $\E$. For linear equations this would mean that the
given equation is trivializable, and the associated curve $L_{\E}$ is
an open part of the rational normal curve in $\RP^n$.

Even though the curve $L_{\E}$ is defined up to contact
transformations, we can define geometric structures on the solution
space itself, which will be independent of the choice of fundamental
solutions. Indeed, let $V(\E)$ be the solution space of the linear
equation~\eqref{lin:eq}. Then at each point $x_0$ on the line we can
define a subspace $V_{x_0}(\E)$ in $V(\E)$ consisting of all solutions
vanishing at $x_0$. Thus, assigning the line $V_{x_0}^{\perp}$ in
the dual space $V(\E)^*$ to each point $x_0$ we get a well-defined
curve in the projective space $\P V_{\E}^*$.

Using the duality principle, we can also construct the curve in the
projective space $\P V_{\E}$. To do this, we can consider the line
$l_{x_0}(\E)$ in $V(\E)$ consisting of all solutions vanishing at
$x_0$ with their derivatives up to order $n-1$. The mapping $x\mapsto
l_{x_0}(\E)$ defines a non-degenerate curve in the projective space
$\P V(\E)$.

Suppose now, that all Wilczynski invariants of the linear
equation~$\E$ vanish identically and consider the symmetry algebra
$\g\subset\sll(V(\E))$ of the constructed curve. Fixing any fundamental
set of solutions and applying Theorem~\ref{wil:thm}, we can show that
this Lie algebra is isomorphic to $\sll(2,\R)$ and acts irreducibly on
$V(\E)$. Let $G$ be the corresponding subgroup in $SL(V(\E)$ and let
$\bar G$ be its product with the central subgroup in $GL(V(\E))$. Then
$\bar G$ is naturally isomorphic to $GL(2,\R)$, and the $\bar G$-module
$V(\E)$ is equivalent to the $GL(2,\R)$-module $S^n(\R^2)$. Thus, any
equation with vanishing Wilczynski invariants defines a
$GL(2,\R)$-structure on its solution space.

For a non-linear equation the vanishing of generalized Wilczynski invariants
is not sufficient for being trivializable. The paper~\cite{d2001}
describes the extra set of invariants that should vanish to guarantee
the trivializability of the equation:
\begin{thm}[\cite{d2001}]
  The equation~\eqref{eq:nonlin} is contact equivalent to the trivial
  equation $y^{(n+1)}=0$ if and only if its generalized Wilczynski
  invariants vanish identically and in addition:
  \begin{itemize}
  \item[] for $n=2$: $f_{2222} = 0$;
  \item[] for $n=3$: $f_{333} = 6f_{233}+f_{33}^2 = 0$;
  \item[] for $n=4$: $f_{44}= 6f_{234}-4f_{333}-3f_{34}^2 =0$;
  \item[] for $n=5$: $f_{55}=f_{45}=0$;
  \item[] for $n\ge 6$: $f_{n,n}=f_{n,n-1}=f_{n-1,n-1}=0$.
  \end{itemize}
\end{thm}

Yet, even if the equation is not trivializable, its linearization at
each solution is trivializable, and we can construct a family of
associated rational normal curves. It will define a
$GL(2,\R)$-structure on the solution space $\S$, if it is a Hausdorff
manifold.

Let $\E$ be an arbitrary non-linear ODE. Suppose that its solution set
$\S$ is Hausdorff (i.e., there exist a quotient of $\E$ by the foliation
formed by all solutions).  Let $y_0(x)$ be any solution of the
non-linear ODE. This is just a point in the manifold $\S$. The
tangent space $T_{y_0}\S$ to $\S$ at $y_0(x)$ can be naturally identified with a
solution space of the linearization of~$\E$ along $y_0(x)$. By the above,
we can naturally construct a curve $l$ in the projectivization of
$T_{y_0}\S$, which means that we have a well-defined two-dimensional
cone $C_{y_0}$ in each tangent space $T_{y_0}\S$.

Suppose now that all generalized Wilczynski invariants of $\E$ vanish
identically. Then linearizations of $\E$ along all solutions are
trivializable, and all cones $C_{y_0}\subset T_{y_0}\S$ are locally equivalent
to the cone in $\R^{n+1}$ corresponding to the normal curve in
$\RP^n$.

Thus, we arrive at the following result.
\begin{thm}\label{thm:gl2}
  Let $\E$ be an arbitrary (non-linear) ODE with vanishing generalized
  Wilczynski invariants. Suppose that its solution space $\S$ is
  Hausdorff and, hence, is a smooth manifold. Then there exists a
  natural irreducible $GL(2,\R)$-structure on $\S$.
\end{thm}

This structure was constructed in~\cite{dun} in a slightly different way and is called a
paraconformal structure.

\subsection{Examples from twistor theory}
In general, explicit examples of equations with vanishing Wilczynski are very
difficult to construct. However, as suggested by N.~Hitchin~\cite{hitchin},
there is a large class of examples coming from twistor theory.

The whole theory above is also valid in complex analytic
category. Consider an arbitrary complex surface $S$ with a rational
curve $L$ on it. Suppose that the normal bundle of $L$ has Grothendieck
type $\O(n)$. Then Kodaira theory states that this rational curve $L$
is included into a complete $(n+1)$-parameter family $\L=\{L_a\}$ of all
deformations of $L$.

The family $\L$ uniquely defines an ordinary differential equation of
order $(n+1)$ such that all these curves are its solutions. In more
detail, we have to consider the lifts of the curves $L_a$ to the jet
space $J^{n+1}(S)$, where they will form a submanifold of codimension
$1$. This defines an ordinary differential equation $\E$.

\begin{thm}
  Let $\E$ be an ordinary differential equation defining the complete
  set of deformations of a rational curve $L$ on a complex surface
  $S$. Then all generalized Wilczynski invariant of $\E$ vanish
  identically.
\end{thm}
\begin{proof}
Each generalized Wilczynski invariant $\W_i$, being restricted to the
lift of any solution $L_a$, defines a global section of
the line bundle $S^iT^*L_a$, which has Grothendieck type $\O(-2i)$. Since
the only global section of the line bundle of this type is zero, we
see that all generalized Wilczynski invariants will vanish identically.
\end{proof}

Consider a number of explicit examples.

\begin{ex}
  Let $S=\CP^2$ and let $L$ be a quadric in $S$. Then the normal
  bundle of $L$ has type $\O(4)$, and the complete family of
  deformations is $5$-dimensional. Clearly, all quadrics in $\CP^2$
  belong to this family and depend on exactly $5$ parameters.

  So, we see that the family $\L$ in this case is a family of all
  quadrics on the complex projective plane. The differential equation
  of all quadrics is well-known and has the form:
\[
  9(y'')^2y^{(5)}-45 y''y'''y^{(4)}+40(y''')^3 = 0.
\]

By above, we have $\W_3=\W_4=\W_5=0$ for this equation. However, this
equation is not trivializable, since, for example, its symmetry
algebra is only $8$-dimensional, while the trivial equation $y^{(5)}=0$
has a $9$-dimensional symmetry algebra.
\end{ex}

\begin{ex}
  Let $S$ be a rational surface $S_k=\mathrm{P}(\O(0)+\O(k))$ viewed
  as a projective bundle over $\CP^1$. Then all sections of this bundle
  with a fixed intersection number $l$ with fibers form a complete
  family of deformations depending on $k+l+1$ parameter. In the
  appropriate coordinate system on $S_k$ they can be explicitly
  written as:
\[
y(x) = \frac{a_0+a_1x+\dots+a_kx^k}{b_0+b_1x+\dots+b_lx^l}.
\]

Assume that $k\ge l$. (Otherwise we can substitute $y$ with $1/y$.)
This set of rational curves defines the following differential
equation:
\[
F_{k,l}=
\begin{vmatrix}
z_{k-l+1} & z_{k-l+2} & \dots & z_{k+1} \\
z_{k-l+2} & z_{k-l+3} & \dots & z_{k+2} \\
\vdots & \vdots & \ddots & \vdots\\
z_{k+1} & z_{k+2} & \dots & z_{k+l+1}
\end{vmatrix} = 0,\qquad\text{where }z_i=\frac{(-1)^i}{i!}y^{(i)}.
\]

For example, in the simplest case $l=0$ we get just the trivial
equation $y^{(k+1)}=0$, while for $l=1$ we get the following equation:
\[
(k+1)y^{(k)}y^{(k+2)}-(k+2)(y^{(k+1)})^2=0.
\]
For $k\ge2$ the symmetry algebra of this equation is
$k+2$-dimensional, and, hence, it is not trivializable.

For $l\ge2$ the equation $F_{k,l}=0$ is never trivializable, since,
for example, being solved with respect to the highest derivative
$y^{(k+l+1}$, it is not linear in term $y^{(k+l)}$.

Since all solutions of the equation $F_{k,l}=0$ are by construction rational
curves, all its generalized Wilczynski invariants vanish identically.
\end{ex}

\begin{ex}
Consider the standard symplectic form $\sigma$ on the 4-dimensional vector
space $\C^4$. Then it defines the contact structure on $\CP^3$
invariant with respect to the induced action of $SP(4,\C)$. It is easy
to see that the set of all rational normal curves in $\CP^3$, which
are at the same time integral curves of the contact distribution,
depends on 7 parameters and forms the complete family of deformations
of any of such curves.

The equation describing this set is given by:
\begin{multline*}
10(y''')^3y^{(7)}-70(y''')^2y^{(4)}y^{(6)}-49(y''')^2(y^{(5)})^2\\
+280y'''(y^{(4)})^2y^{(5)}-175(y^{(4)})^4=0.
\end{multline*}

Again, since all solutions of this equation are rational curves, all
its generalized Wilczynski invariants vanish identically. Yet, this
equation is not trivializable, since its symmetry algebra is
10-dimensional (in fact, it coincides with $\mathfrak{sp}(4,\C)$),
while the trivial equation $y^{(7)}=0$ has 11-dimensional symmetry
algebra.
\end{ex}

Note that equations from the examples above have the common property that they have the symmetry
algebra of submaximal dimension~\cite{olver}.

\section{Wilczynski invariants for systems of ODEs}
\label{sec:sys}
Since all results of this paper are based on two ideas, namely the
notion of linearization and the invariants and structures associated
with linear equations, they are directly generalized to systems of
ordinary differential equations and even to equations of finite type.

The generalization of Wilczynski invariants to the systems of ODEs
was obtained by Se-ashi~\cite{se2}. Consider an arbitrary system of
linear ordinary differential equations:
\[
y^{(n+1)}+P_n(x)y^{(n)}+\dots+P_0(x)y(x)=0,
\]
where $y(x)$ is an $\R^m$-valued vector function. The canonical
Laguerre-Forsyth form of these equations is defined by conditions
$P_{n}=0$ and $\tr P_{n-1}=0$. Then, as in the case of a singe ODE,
the following expressions:
\begin{equation}
\Theta_k = \sum_{j=1}^{k-1} (-1)^{j+1}
\frac{(2k-j-1)!(n-k+j)!}{(k-j)!(j-1)!} P_{n-k+j}^{(j-1)}, \quad
k=2,\dots, n+1.
\end{equation}
are the $\End(\R^m)$-valued relative invariants, where each invariant
$\th_i$ has weight $i$. Note that unlike the case of a single ODE, the
first non-trivial Wilczynski invariant has weight $2$.

We can also consider a characteristic curve associated with any linear
system, which will take values in the Grassmann manifold
$\Gr_{m}(\R^{(n+1)m})$.  Se-ashi also proved that the system is
trivializable if and only if all these invariants vanish identically,
or, equivalently, if the characteristic curve is an open part of the
curve defined by a trivial equation. The symmetry group of this curve
is isomorphic to the direct product of $SL(2,\R)$ and $GL(m,\R)$ with
an action on $\R^{(n+1)m}$ equivalent to the natural action on
$S^{(n+1)}(\R^2)\otimes \R^m$.

All that leads us immediately to the following generalization of
Theorems~\ref{thm:winv} and~\ref{thm:gl2}.
\begin{thm}\label{thm:sysode}
  Let $\E$ be an arbitrary (non-linear) system of $m$ ordinary
  differential equations of order $(n+1)$. Then there exist generalized
  Wilczynski invariants $\W_2$, \dots, $\W_{n+1}$ that, being
  restricted to each solution, coincide with Wilczynski invariants of
  the linearization along this solution.

  Suppose that all these invariants vanish identically and the
  solution space $\S$ is a smooth manifold. Then there exists a natural
  irreducible $SL(2,\R)\times GL(m,\R)$ structure on~$\S$.
\end{thm}

\begin{rem}
Let us note that in the case of systems of ODEs the invariant $\W_i$ is a
section the vector bundle $S^iE^*\otimes\End(V_n)$, where $V_n$ is a
subbundle of the normal bundle $N(\E)=T\E/E$ defined as
$\pi_{n+1,n}^*C^n\cap T\E$ modulo $E$. In case of a single ODE the bundle $V_n$
is one-dimensional, and $\End(V_n)$ can be identified with the trivial
bundle over $\E$.
\end{rem}

\begin{ex}
  Consider the simplest non-trivial case of second order systems of
  ODEs. The linear homogeneous system can be written as:
  \begin{equation}\label{linsyst}
  y''=A(x)y'+B(x)y,
  \end{equation}
  where $y(x)$ is the unknown $\R^m$-valued function and $A(x),B(x)\in
  \End(\R^m)$. The set of all transformations preserving the class of
  such systems has the form $(x,y)\mapsto (\lambda(x),\mu(x)y)$, where
  $\lambda(x)$ is a local line reparametrization and $\mu(x)\in
  GL(m,\R)$.

  We can bring the equation~\eqref{linsyst} into the semi-canonical
  form with vanishing coefficient $A(x)$ by means of a certain
  gauge transformation $(x,y)\mapsto (x,\mu(x)y)$. Then, applying a
  reparametrization $(x,y)\mapsto(\lambda(x),y)$ and an
  appropriate gauge transformation to bring the equation back to the
  semi-canonical form, we can make trace of $B(x)$ vanish.

  The traceless part of the coefficient $B(x)$ in the semicanonical
  form gives us the Wilczynski invariant $\Th_2$ for linear systems of
  second order. Explicitly, it can be written as:
  \begin{equation}\label{wi2ord}
  \Th_2 = \Phi(x)-1/m\tr \Phi(x),\quad\text{where
  }\Phi(x)=B(x)-\frac{1}{2}A'(x)+\frac{1}{4}A(x)^2.
  \end{equation}
  According to Se-ashi~\cite{se2}, this is the only Wilczynski
  invariant available for systems of second order, and the
  equation~\eqref{linsyst} is trivializable if and only if the
  invariant $\Th_2$ vanishes identically.

  Consider now an arbitrary non-linear system of second order:
  \[
  y_i''=f_i(x,y_j,y_k'),\qquad i=1,\dots,m.
  \]
  Then its generalized Wilczynski invariant can be obtained
  from~\eqref{wi2ord} by substituting $A(x)$ with the matrix
  $\left(\frac{\partial f_i}{\partial y_k'}\right)$, $B(x)$ with the
  matrix $\left(\frac{\partial f_i}{\partial y_j}\right)$, and the
  usual derivative $d/dx$ with operator of total derivative:
  \[
  D=\pp{x}+\sum_{i=1}^m y_i'\pp{y_i}+\sum_{j=1}^m f_j\pp{y_j'}.
  \]
  Denote by $W_2$ the invariant we get in this way. It coincides (up
  to the constant multiplier) with the invariant constructed by M.~Fels~\cite{fels}.

  If the invariant $W_2$ vanishes identically, then we get the
  $GL(2,\R)\times GL(m,R)$ structure on the solution space. This
  structure is also called Segre structure and was first constructed for
  solution space of second order ODEs with vanishing generalized
  Wilczynski invariant by D.~Grossman~\cite{gros}.
\end{ex}

\end{document}